\documentclass[11pt]{article}
\usepackage{amssymb,amsmath,mathrsfs,txfonts,graphicx,tikz,color}
\graphicspath{{figs/}}

\begin{document}

\allowdisplaybreaks
\def\pd#1#2{\frac{\partial#1}{\partial#2}}
\def\dfrac{\displaystyle\frac}
\let\oldsection\section
\renewcommand\section{\setcounter{equation}{0}\oldsection}
\renewcommand\thesection{\arabic{section}}
\renewcommand\theequation{\thesection.\arabic{equation}}

\newtheorem{theorem}{\indent Theorem}[section]
\newtheorem{lemma}{\indent Lemma}[section]
\newtheorem{proposition}{\indent Proposition}[section]
\newtheorem{definition}{\indent Definition}[section]
\newtheorem{remark}{\indent Remark}[section]
\newtheorem{corollary}{\indent Corollary}[section]

\title{\LARGE
Variational Approach of  Critical Sharp Front Speeds in
Density-dependent Diffusion Model with Time Delay
}
\author{
Tianyuan Xu$^{a,d}$, Shanming Ji$^{b,d}$\thanks{Corresponding author, email:jism@scut.edu.cn},
Ming Mei$^{c,d}$, Jingxue Yin$^a$,
\\
\\
{ \small \it $^a$School of Mathematical Sciences, South China Normal University}
\\
{ \small \it Guangzhou, Guangdong, 510631, P.~R.~China}
\\
{ \small \it $^b$School of Mathematics, South China University of Technology}
\\
{ \small \it Guangzhou, Guangdong, 510641, P.~R.~China}
\\
{ \small \it $^c$Department of Mathematics, Champlain College Saint-Lambert}
\\
{ \small \it Quebec,  J4P 3P2, Canada, and}
\\
{ \small \it $^d$Department of Mathematics and Statistics, McGill University}
\\
{ \small \it Montreal, Quebec,   H3A 2K6, Canada}
}
\date{}

\maketitle

\begin{abstract}
For the classical reaction diffusion equation,
the priori speed of fronts
is determined exactly in the pioneering paper
(R.D. Benguria and M.C. Depassier, {\em Commun. Math. Phys.} 175:221--227, 1996)
by variational characterization
method.
In this paper, we model the dispersal process using a density-dependent
diffusion equation with time delay.
We show the existence and uniqueness of sharp critical fronts,
where the sharp critical front is $C^1$-smooth when
the diffusion degeneracy is weaker with $1<m<2$,
and the sharp critical front is non-$C^1$-smooth (piecewise smooth)
when the diffusion degeneracy is stronger with $m\ge 2$.
We give a new variational approach for the critical wave speed and
investigate how the time delay affects the
propagation mechanism of fronts.
Our results provide some interesting insight into the
dynamics of critical traveling wave.


\end{abstract}


\section{Introduction and Preliminaries}
We are interested in the critical sharp traveling wave for the following density-dependent
diffusion equation with time delay
\begin{align}\label{eq-main}
\pd u t=D\Delta u^m-d(u)+b(u(t-r,x)),\quad x\in \mathbb R,~t>0,
\end{align}
where $u$ is the population density, $b(u(t-r,x))$ is the birth function,
$r\ge0$ is the time delay, $D>0$ represents  the diffusivity of the population,
and $d(u)$ is the death rate function.
The diffusion term $D\Delta u^m$~(with $m>1$) is considered to be in the form of porous medium type,
which is dependent on the population density
due to the population pressure \cite{Gurney75,Murry,Shiguesada79}.

Our main purpose of this paper is
to study the critical sharp traveling waves for \eqref{eq-main}.
A traveling wave solution is a special solution in the
form of $u(t, x)=\phi(x+ct)$ with wave speed  $c$.
We show the existence and uniqueness of critical traveling
wave for \eqref{eq-main}, and investigate the geometric shape of this wave, which is
partially compactly supported (called ``sharp wave'').
Inspired by the work of Benguria and Depassier \cite{Benguria} for the regular diffusion equation,
we give a new variational characterization of this sharp critical wave for the density-dependent
degenerate diffusion equation with time delay.
We demonstrate that, the diffusion degeneracy and the time delay
both act as two important mechanisms for the propagation of wave front.
The speed of this sharp critical front is reduced due to the time delay effect.
The critical wave speed $c^*(m,D,r)$ (defined in \eqref{eq-def})
of model \eqref{eq-main} with time delay is smaller than
the critical wave speed $c^*(m,D,0)$ of the case without time delay as shown
in Theorem \ref{th-existence},
and when the diffusion degeneracy is weaker with $1<m<2$, the
wave front is still $C^1$-smooth,
while, when the diffusion degeneracy is stronger with $m\ge2$,
then the wave front is non-$C^1$-smooth.


For functions $d(s)$ and $b(s)$, we have the following hypotheses:
\begin{align} \nonumber
&\text{There exist~}u_-=0, u_+>0 \text{~such that~}
d, b\in C^2([0,u_+]),
d(0)=b(0)=0,
\\ \label{eq-H}
&d(u_+)=b(u_+),
b'(0)>d'(0)\ge 0, d'(u_+)\ge b'(u_+)\ge0,
d'(s)\ge 0, b'(s)\ge0.
\end{align}

Here,
both $u_-=0$ and $u_+>0$
are constant equilibria of \eqref{eq-main},
and functions $b(u)$,
$d(u)$ are both nondecreasing.
The assumption \eqref{eq-H} is summarized from the classical Fisher-KPP
equation \cite{Fisher}, see also a lots of evolution equations in ecology,
for example,
the well-studied Nicholson's blowflies equation \cite{Gurney80} with
the death function
$d_1(u)=\delta u$ or $d_2(u)=\delta u^2$,
the birth function
$$b_1(u)=pue^{-au^q},
\quad p>0,\quad q>0,\quad a>0;$$
and the Mackey-Glass equation \cite{Mackey} with
the growth function
$$
b_2(u)=\frac{pu}{1+au^q}\quad p>0,\quad q>0,\quad a>0.
$$

Since the pioneering work of Schaaf \cite{Schaaf} on traveling fronts
in delayed reaction-diffusion equations,
this field has been extensively investigated \cite{Alomari,Gomez,So}.
So and Zou \cite{So} studied travelling front
solutions of the Nicholson's blowflies equation.
They showed that a monotone wave front exists connecting the
equilibria for monotone birth function.
In \cite{Faria} the existence of non-monotone travelling
fronts of delayed reaction-diffusion equation is proved,
which includes several specific classes of birth rate functions.
For the stability of critical and noncritical traveling waves,
we refer the readers to \cite{Chern,MeiSIAM14,Mei_LinJDE09,MeiJDE09,MeiSIAM10,IJNAM}
and the references therein.

Density-dependent dispersal has been observed in many biological populations \cite{Gurney75,Murry,Shiguesada79}.
Individuals migrate from densely populated areas to sparsely areas
to avoid overcrowding \cite{Okubo}.
The probability of leaving the current site increases with the local population density.
This density-dependent dispersal mechanism arises from
competition between conspecifics or deteriorating environmental conditions \cite{Matthysen}.

One of the most interesting features of the systems with density dependent diffusion is the existence of sharp type traveling wave.
Sharp type traveling waves are zero on a half-plane and decay to zero in a continuous but non-smooth way. The interest in sharp traveling waves is related to the important property of finite speed of propagation, as showed in \cite{PME}.
The sharp traveling wave
plays an important role
in the analysis of the propagation properties of degenerate diffusion equations
since the equation generates compactly supported solutions if the
initial value is compactly supported.
The appearance of sharp profiles was first discussed in the
pioneering work of Aronson \cite{Aronson80}, which is within the framework of models for density dependent biological invasion.
In the case without time delay, the sharp type traveling wave
is unique, and the corresponding speed is the minimal admissible traveling wave speed
of all types, and further the speed is also the spreading speed,
see the very recent paper by Audrito and V\'azquez \cite{VazquezJDE17}
for a doubly nonlinear diffusion equation.

When the system is with a time time delay $r>0$, the situation changes dramatically.
Time delay and degenerate diffusion
lead to essential difficulties in proving the similar results as in \cite{VazquezJDE17} concerning
the existence, uniqueness and critical
properties of sharp traveling waves for delayed system.
To the best of our knowledge, the study of traveling waves for the monostable delayed model \eqref{eq-main} with density dependent diffusion was initiated in \cite{HJMY,JDE18XU}.
But the existence of the critical traveling waves $\phi(x+c^*t)$ still remains open.
It is proved in \cite{JDE18XU}
that the admissible traveling wave speeds for \eqref{eq-tw} are greater than
or equal to $c^*(m,D,r)>0$ if the time delay is small
(see Theorem 2.4 in \cite{JDE18XU}),
but nothing is known about the type of the critical traveling waves.
If the birth rate function is not restricted to be monotone,
the authors in \cite{HJMY} proved that \eqref{eq-tw} admits
traveling waves for some $c(r)$ if the time delay is small
but those speeds are not critical (see Theorem 1.2 in \cite{HJMY}).
In our recent work \cite{non-monotone}, we find the existence of
sharp type traveling wave (wavefront or semi-wavefront).
However, the uniqueness of such waves and
whether it is corresponding to the minimal admissible wave speed still remain open.
In this paper, we will answer those unsettled questions
in \cite{HJMY,JDE18XU,non-monotone} under the hypotheses \eqref{eq-H}.

Our main results  provide some interesting insight into
the dynamics of critical traveling wave solutions for \eqref{eq-main}.
Theorem \ref{th-existence} stated below shows that
model \eqref{eq-main} admits a unique sharp traveling wave and the
sharp traveling wave is monotonically increasing.
Here, we develop a phase transform approach to show
more precise behavior about traveling waves.
Generally speaking, since the trajectories may intersect with each other due to the time delay,
traditional phase plane analysis method is incapable of showing the existence of solutions
of time delayed model.
Surprisingly, we find that it provides a blueprint to draw more precise information
about the solutions.

Critical traveling wave plays an important role in biological invasion.
For this degenerate diffusion equation with time delay,
the critical traveling wave is of sharp type.
This sharp traveling wave $\phi(x+c^*t)$ separates the existence or the nonexistence of wavefronts.
When $c<c^*$, there is no traveling wave for \eqref{eq-main}.
When $c>c^*$, all the traveling waves are smooth and positive.

It is worthy of mentioning that the time delay leads to
speed reducing mechanism of critical wave speed.
For the case without time delay and with linear diffusion (i.e. $m=1$ and $r=0$),
it is proved by Benguria and Depassier \cite{Benguria} that
$$
c^*(1,D,0)=\max\left\{2\sqrt{D(b'(0)-d'(0))},
\sup_{g\in \mathscr{D}}J_1(g)\right\},
$$
where $J_1$ and $\mathscr{D}$ are defined in \eqref{eq-cstar0}.
As far as we know, this is the first literature that provides an effective
method to calculate the velocity of the fronts.
A well-known result for the classical Fisher-KPP equation $u_t-D\Delta u=u(1-u/K)$
shows that $c^*(1,D,0)=2\sqrt{D}$ for this special types of
birth and death functions.
In this paper,
we give a new variational characterization of the critical wave speed
for this density dependent model with time delay
and show the effect of time delay to the
critical wave speed.
The critical wave speed $c^*(m,D,r)$ of model \eqref{eq-main} with time delay is smaller than
the critical wave speed $c^*(m,D,0)$ of the case without time delay as
shown in Theorem \ref{th-existence}.

The rest of this paper is organized as follows.
In section 2, we present the main results on the existence, uniqueness and critical
properties of sharp traveling wave solutions.
Section 3 is devoted to the proof of the existence and properties of critical traveling
wave solutions. We show that the critical wave for model \eqref{eq-main} is sharp type and it
is unique. Then we prove the effect of time delay to the critical wave speed and
the regularity of the sharp traveling wave.

\section{Main results}
We consider the degenerate diffusion equation with time delay \eqref{eq-main}.
We are looking for the traveling wave solutions that
connect the two equilibria $u_-=0$ and $u_+=:K$.
Under the hypotheses \eqref{eq-H},
the birth function $b(u)$ is monotonically increasing on $[u_-,u_+]=[0,K]$.
Let $\phi(\xi)$, where $\xi=x+ct$ and $c>0$,
be the traveling wave solution of  \eqref{eq-main},
we get (we write $\xi$ as $t$ for the sake of simplicity)
\begin{align}\label{eq-tw}
\begin{cases}
\displaystyle
c\phi'(t)=D({\phi^m}(t))''-d(\phi(t))+b(\phi(t-cr)),\quad t\in\mathbb R,\\
\phi(-\infty)=0, \quad \phi(+\infty)=K.
\end{cases}
\end{align}

The traveling wave with partially compact support (we will call this ``sharp type'')
plays a crucial role
in the analysis of the propagation properties of degenerate diffusion equations.
In the case without time delay, the sharp type traveling wave
is unique, and the corresponding speed is the minimal admissible traveling wave speed
of all types, and further the speed is also the propagation speed,
see the interesting paper by Audrito and V\'azquez \cite{VazquezJDE17}
for a doubly nonlinear diffusion equation
and the pioneering work by Aronson \cite{Aronson80} for the
porous medium equation.
Here, we try to show the uniqueness of sharp type traveling wave
and the critical property of the corresponding wave speed for the time delayed case.
And further we will prove that the time delay slows down the
critical wave speed.

\begin{figure}[htb]
\begin{center}
\begin{tikzpicture}[scale=0.8,domain=-6:6]
\def\axl{-6} \def\axr{6} \def\ay{4} \def\pa{\ay/180} \def\pb{3} \def\pc{2}
\def\qa{-0.1531} \def\qb{1.1531} \def\ra{1.3606} \def\rb{1.2789}
\draw[->,>=latex,line width=.5pt] ({\axl},0)--(\axr,0) node[right] {$\xi$};
\draw[->,>=latex,line width=.5pt] (0,-0.2)--(0,{\ay+0.3}) node[right] {$\phi$};
\node[below] at (-.2,.05) {$O$};
\draw[dashed, line width=.6pt] ({\axl},\ay)--(\axr,\ay);
\draw[color=blue,line width=.7pt]
    plot[domain=-6:6, samples=144, smooth] (\x,{\ay/2+\pa*atan(\pb*\x)});
\node at (0.9,3.2) {($a$)};
\draw[color=red,line width=.7pt]
    plot[domain=2:6, samples=144, smooth] (\x,{\ay/2+\pa*atan(\pb*(\x-2))});
\draw[color=red,line width=.7pt]
    plot[domain=0:1, samples=144, smooth] (\x,{\x*\x/3});
\draw[color=red,line width=.7pt]
    plot[domain=1:2, samples=144, smooth]
    (\x,{1/3+2/3*(\x-1)+\qa*(\x-1)*(\x-1)+\qb*(\x-1)^3});
\draw[color=red,line width=.7pt]
    plot[domain=1:6, samples=144, smooth] (\x,{\ay/2+\pa*atan(\pb*(\x-1))});
\draw[color=red,line width=.7pt]
    plot[domain=0:0.5, samples=144, smooth] (\x,{\x*3/2});
\draw[color=red,line width=.7pt]
    plot[domain=0.5:1, samples=144, smooth]
    (\x,{3/4+3/2*(\x-1/2)+\ra*(\x-1/2)*(\x-1/2)+\rb*(\x-1/2)^3});
\draw[color=red,line width=.7pt] ({\axl},0)--(0,0);
\node at (1.3,2.4) {($b_1$)};
\node at (2.2,1.8) {($b_2$)};
\end{tikzpicture}
\end{center}
\caption{Traveling waves: ($a$) smooth type;
($b_1$) non-$C^1$ sharp type; ($b_2$) $C^1$ sharp type.}
\label{fig-tw}
\end{figure}
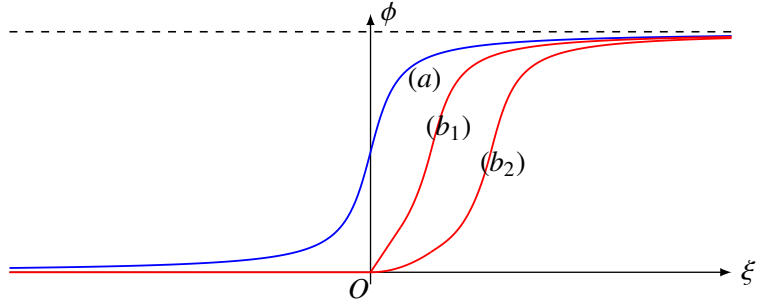

We are looking for traveling waves of sharp type
and try to show the relation between the corresponding wave speed
and all the admissible traveling wave speeds.
We present the following definition of
sharp and smooth traveling waves, see Figure \ref{fig-tw} for illustration.
Here are some notations used throughout the paper:
$L_\text{loc}^1(\mathbb R)$ is the set of locally Lebesgue integrable functions,
$C^1=C^1(\overline{\mathbb R})$,
$$
C_\text{unif}^b(\mathbb R):=
\{\phi\in C(\mathbb R)\cap L^\infty(\mathbb R);\phi \text{~is uniformly continuous on~}\mathbb R\},
$$
and
$$W_\text{loc}^{1,2}(\mathbb R):=\{\phi; \phi\in W^{1,2}(\Omega)
\text{~for any compact subset~}\Omega\subset\mathbb R\}.$$

\begin{definition} \label{de-semi}
A profile function $\phi(t)$ is said to be a traveling wave solution (TW)
of \eqref{eq-tw} if $\phi\in C_\text{unif}^b(\mathbb R)$, $0\le \phi(t)\le K:=u_+$,
$\phi(-\infty)=0$, $\phi(+\infty)=K$,
$\phi^m\in W_{\text{loc}}^{1,2}(\mathbb R)$,
$\phi(t)$ satisfies \eqref{eq-tw} in the sense of distributions.
The TW $\phi(t)$ is said to be of sharp type if
the support of $\phi(t)$ is semi-compact, i.e.,
$\text{supp}\,\phi=[t_0,+\infty)$ for some $t_0\in\mathbb R$,
$\phi(t)>0$ for $t>t_0$.
On the contrary, the TW $\phi(t)$
is said to be of smooth type if $\phi(t)>0$ for all $t\in\mathbb R$.

Furthermore, for the sharp TW $\phi(t)$,
if $\phi''\not\in L_\text{loc}^1(\mathbb R)$,
we say that $\phi(t)$ is a non-$C^1$ type sharp TW;
otherwise, if $\phi''\in L_\text{loc}^1(\mathbb R)$,
we say that $\phi(t)$ is a $C^1$ type sharp TW.
\end{definition}

Without loss of generality, we may always shift $t_0$ to $0$
for the sharp type TW.
Therefore, a sharp type TW $\phi(t)$ is
a special TW such that $\phi(t)=0$ for $t\le0$, and $\phi(t)>0$ for $t>0$.

For any given $m>1$, $D>0$, and $r\ge0$, we define the critical (or minimal) wave speed
$c^*(m,D,r)$ for the degenerate diffusion equation \eqref{eq-tw} as follows
\begin{equation} \label{eq-def}
c^*(m,D,r)
:=\inf\{c>0; \eqref{eq-tw} \text{~admits increasing TWs with speed $c$}\}.
\end{equation}
For the case without time delay and with degenerate diffusion (i.e. $m>1$ and $r=0$),
it is proved in \cite{HJMY} that
\begin{equation} \label{eq-cstar0}
c^*(m,D,0)=\sup_{g\in \mathscr{D}}\left\{2\sqrt{D}\int_0^K
\sqrt{-ms^{m-1}g(s)g'(s)(b(s)-d(s))}ds\right\}
=:\sup_{g\in \mathscr D}J_m(g),
\end{equation}
where $\mathscr{D}=\{g\in C^1([0,K]);g(K)=0,\int_0^K g(s)ds=1,g'(s)<0,\forall s\in(0,K)\}$.
We note that in \cite{HJMY}, we prove the above variational characterization
for the Nicholson's blowflies model and we can verify that it holds true for the general type.
For the case without time delay and with linear diffusion (i.e. $m=1$ and $r=0$),
it is proved by Benguria and Depassier \cite{Benguria} that
$$
c^*(1,D,0)=\max\left\{2\sqrt{D(b'(0)-d'(0))},
\sup_{g\in \mathscr{D}}J_1(g)\right\}.
$$
A well-known result for the classical Fisher-KPP equation $u_t-D\Delta u=u(1-u/K)$
shows that $c^*(1,D,0)=2\sqrt{D}$ for this special type of
birth and death functions.

In this paper we show that \eqref{eq-tw} admits a unique
sharp type TW, and the sharp TW is monotonically increasing
and corresponding to the minimal wave speed $c^*(m,D,r)$,
and further $c^*(m,D,r)<c^*(m,D,0)$ for any time delay $r>0$.
As a consequence, the time delay slows down the minimal traveling wave speed
for the degenerate diffusion case.

Our main results are as follows.

\begin{theorem}[Critical Sharp Traveling Wave] \label{th-existence}
Assume that $d(s)$ and $b(s)$ satisfy \eqref{eq-H}, and $m>1$, $r\ge0$.
There exists a unique $c^*=c^*(m,D,r)>0$ defined in \eqref{eq-def}
satisfying $c^*(m,D,r)<c^*(m,D,0)$ for any time delay $r>0$,
such that \eqref{eq-tw} admits a unique (up to shift)
sharp traveling wave $\phi(x+c^*t)$ with speed $c^*$,
which is the critical traveling wave of \eqref{eq-tw} and is monotonically increasing.
\end{theorem}

\begin{theorem}[Smooth Traveling Waves] \label{th-smooth}
Assume that the conditions in Theorem \ref{th-existence} hold.
For any $c>c^*=c^*(m,D,r)>0$ defined in \eqref{eq-def},
the traveling wave $\phi(x+ct)$ of \eqref{eq-tw} with speed $c$ is smooth type
for all $m>1$ such that $\phi\in C^2(\mathbb R)$;
while for any $c\le c^*=c^*(m,D,r)$,
there is no smooth type traveling waves of \eqref{eq-tw} with speed $c$.
\end{theorem}

\begin{remark}
Theorem \ref{th-existence} and Theorem \ref{th-smooth} show that
the sharp traveling wave $\phi(x+c^*t)$ is unique and critical in the following senses:
\\ \indent
(i) the speed $c^*$ is the unique wave speed such that \eqref{eq-tw}
admits sharp type traveling waves (no other wave speed for sharp traveling waves);
\\ \indent
(ii) the speed $c^*$ is critical such that no traveling waves exist for $c<c^*$
and $\phi(x+c^*t)$ is the wave with minimal wave speed;
\\ \indent
(iii) $\phi(x+c^*t)$ is the unique traveling wave with speed $c^*$,
i.e., no smooth traveling waves with the critical wave speed exist.
\\ \indent
Therefore, a traveling wave of \eqref{eq-tw} is sharp if and only if the
corresponding wave speed is critical.
\end{remark}

\begin{remark}
The existence of non-critical traveling waves was proved in \cite{JDE18XU},
where no regularity results were given as the types of traveling waves,
sharp or smooth, were not settled.
We also note that the $C^2$ regularity and existence of traveling waves
were shown in \cite{HJMY} for some wave speeds $c>c^*(m,D,0)$,
which according to our results here (Theorem \ref{th-existence})
are non-critical since the critical wave speed $c^*(m,D,r)<c^*(m,D,0)$.
Theorem \ref{th-smooth} presents the $C^2$ regularity of non-critical traveling waves
with any speed $c>c^*(m,D,r)$.
\end{remark}

The sharp traveling wave is classified into
$C^1$ type and non-$C^1$ type according to the degeneracy index $m$.

\begin{theorem}[Regularity of Sharp Wave] \label{th-sharp}
Assume that the conditions in Theorem \ref{th-existence} hold.
If $m\ge2$, then the sharp traveling wave is of non-$C^1$ type;
while if $1<m<2$, then the sharp traveling wave is of $C^1$ type.
\end{theorem}

\begin{remark}
Regarding the regularity of sharp waves,
roughly speaking the degeneracy strengthens as $m>1$ increases
and the regularity of the case $m\ge2$ is weaker than that of $1<m<2$.
For the case $1<m<2$, the sharp traveling wave remains $C^1$ regularity
but not analytic.
\end{remark}

\section{Proof of the main results}

In this section, we first solve \eqref{eq-tw} locally for sharp type solutions,
and then we develop a phase transform approach
to show the monotone dependence of the wave speed.
The sharp traveling wave is the unique local solution that exists globally,
monotonically increases on the whole real line, and has the least upper bound $K$.
The corresponding wave speed is characterized via a variational approach
inspired by Benguria and Depassier \cite{Benguria}
(see also Huang et al. \cite{HJMY}),
and is also compared with the smooth type traveling wave speeds.

The existence of sharp TW for the case without time delay
is proved in \cite{HJMY} for the Nicholson's blowfly model.
It is also valid for the general birth rate and death rate functions without time delay
and here we only focus on the case with time delay $r>0$.

For any given $D>0$ and $r>0$,
we solve \eqref{eq-tw} step by step.
First, noticing that the sharp wave solution $\phi(t)=0$ for $t\le0$
and then $\phi(t-cr)=0$ for $t\in[0,cr)$, \eqref{eq-tw} is locally reduced to
\begin{equation} \label{eq-semi-1}
\begin{cases}
c\phi'(t)=D(\phi^m(t))''-d(\phi(t)), \\
\phi(0)=0, \quad (\phi^m)'(0)=0, \quad t\in(0,cr),
\end{cases}
\end{equation}
whose solutions are not unique and we choose the maximal one
such that $\phi(t)>0$ for $t\in(0,cr)$ as shown in the following lemma.
Here, $(\phi^m)'(0)=0$ is a necessary and sufficient condition
such that the zero extension of $\phi(t)$ to the left
satisfies \eqref{eq-tw} locally near $0$ in the sense of distributions.

The following three lemmas are formulated and proved in our related paper
\cite{non-monotone} in order to show the existence of sharp TWs
in a more general setting.
We omit the proofs here for the sake of simplicity.

\begin{lemma}[\cite{non-monotone}] \label{le-semi-1}
For any $c>0$, the degenerate ODE \eqref{eq-semi-1} admits a unique maximal solution
$\phi_c^1(t)$ on $(0,cr)$ such that $\phi_c^1(t)>0$ on $(0,cr)$ and
$$
\phi_c^1(t)=\Big(\frac{(m-1)c}{Dm}t\Big)^\frac{1}{m-1}+o(t^\frac{1}{m-1}),
\quad t\to0^+.
$$
\end{lemma}

Next, let $\phi_c^2(t)$ be the solution of the following initial value ODE problem
\begin{equation} \label{eq-semi-2}
\begin{cases}
c\phi'(t)=D(\phi^m(t))''-d(\phi(t))+b(\phi_c^1(t-cr)), \\
\phi(r)=\phi_c^1(r), \quad
\phi'(r)=(\phi_c^1)'(r), \qquad t\in(cr,2cr).
\end{cases}
\end{equation}
The above steps can be continued unless $\phi_c^k(t)$ blows up
or decays to zero in finite time for some $k\in \mathbb N^+$.
Let $\phi_c(t)$ be the connecting function of those functions on each step, i.e.,
\begin{equation} \label{eq-semi}
\phi_c(t)=
\begin{cases}
\phi_c^1(t), \quad &t\in[0,cr),\\
\phi_c^2(t), \quad &t\in[cr,2cr),\\
\dots\\
\phi_c^k(t), \quad &t\in[(k-1)cr,kcr),\\
\dots
\end{cases}
\end{equation}
for some finite steps such that $\phi_c(t)$ blows up or decays to zero,
or for infinite steps such that $\phi_c(t)$ is defined on $(0,+\infty)$
and zero extended to $(-\infty,0)$ for convenience.

\begin{lemma}[\cite{non-monotone}] \label{le-semi-decay}
For any given $m$, $D$ and $r>0$, there exists a constant $\underline c>0$ such that
if $c\le \underline c$, then
$\phi_c(t)$ decays to zero in finite time.
\end{lemma}

\begin{lemma}[\cite{non-monotone}] \label{le-semi-blow}
For any given $m$, $D$ and $r>0$, there exists a constant $\overline c>0$ such that
if $c\ge \overline c$, then
$\phi_c(t)$ grows up to $+\infty$ as $t$ tends to $+\infty$.
\end{lemma}

In the paper \cite{non-monotone}, the birth rate function is not restricted to be monotone
in $[0,K]$, and the time delay together with the non-monotone structure
of birth rate function gains the possibility of the existence of
oscillating sharp TWs.
As a result, the monotone dependence of $\phi_c(t)$ with respect to $c$
is generally not true, and the uniqueness of the
wave speed for sharp type TWs remains open under these general settings.

Here in this paper, under the assumption \eqref{eq-H},
the birth rate function is monotonically increasing in $[0,K]$
and we can present a positive answer to the above questions.

\begin{lemma}[Continuous Dependence] \label{le-continu}
For any given $m$, $D$ and $r>0$,
the solution $\phi_c(t)$ is locally continuously dependent on $c$.
That is, for any $c>0$ and any given $T>0$ and $\varepsilon>0$,
there exists a $\delta>0$ such that for any $|c_1-c|<\delta$ and $c_1>0$ we have
$$
|\phi_{c_1}(t)-\phi_c(t)|<\varepsilon,
\quad \forall t\in(0,T_1-\varepsilon),
$$
where $T_1=\min\{T,T_c\}$
with $T_c$ being the existence interval of $\phi_c(t)$.
\end{lemma}
{\it\bfseries Proof.}
Without loss of generality, we may assume that $T_1>cr+\varepsilon$.
The proof is divided into two parts:
the continuous dependence of the singular ODE \eqref{eq-semi-1} within $(0,cr)$
and the continuous dependence of a regular ODE within $(cr,T_1-\varepsilon)$.

Step I. We prove that $\phi_c(t)$ together with $\phi_c'(cr)$ is
continuously dependent on $c$ for $t\in(0,cr)$.
Since $T_c>cr$, we see that $\phi_c(t)$ is positive for $t\in(0,cr]$.
We note that the maximal solution $\phi_c(t)$ is the unique solution such that
$\phi_c(t)>0$ in a right neighbor of $0$
and the asymptotic analysis Lemma \ref{le-semi-1} shows that
$\phi_c(t)$ is locally monotonically and continuously dependent on $c$
within some interval $(0,t_1)\subset (0,cr)$.
In $(t_1,cr)$, $\phi_c(t)$ is bounded away from zero
and \eqref{eq-semi-1} is a regular ODE,
and the continuously dependence follows from the classical theory.

Step II. We prove that $\phi_c(t)$ is
continuously dependent on $c$ for $t\in(cr,T_1-\varepsilon)$.
As $T_c$ is the existence interval of $\phi_c(t)$,
$\phi_c(t)$ is bounded from above and below
for $t\in[cr,T_c-\varepsilon]$ such that
$\phi_c(t)\in[M_1,M_2]\subset(0,+\infty)$ for some $M_2>M_1>0$.
According to the construction of $\phi_c(t)$, we see that $\phi_c(t)$ satisfies
\begin{equation*}
\begin{cases}
c\phi'(t)=D(\phi^m(t))''-d(\phi(t))+b(\phi(t-cr)), \quad t\in(cr,T_1-\varepsilon),\\
\phi(cr)=\phi_c(cr), \quad \phi'(cr)=\phi_c'(cr),
\end{cases}
\end{equation*}
which is a regular ODE without singularity on bounded interval.
This completes the proof.
$\hfill\Box$

\begin{lemma} \label{le-monotone}
For any given $c>0$, $\phi_c(t)$ is locally strictly increasing
on $(0,t_*)$ for some $t_*\in (0,+\infty]$.
We may slightly abuse the notation and
denote the maximal interval in which $\phi_c(t)$ is strictly increasing
and $\phi_c(t)<K$ also by $(0,t_*)$.
If $t_*<+\infty$ and $\phi_c(t_*)<K$,
then $\phi_c(t)$ decreases to zero
and will never grows up to $K$ after $t_*$.
\end{lemma}
{\it\bfseries Proof.}
The local monotonicity follows from Lemma \ref{le-semi-1}.
We first show that $\phi_c(t)$ can not be constant after $t_*$.
Otherwise,
$$
b(\phi(t-cr))=c\phi'(t)-D({\phi^m}(t))''+d(\phi(t))=d(\phi(t))
$$
is a constant after $t_*$, which contradicts
to the strictly monotone increasing of $\phi_c(t)$ before $t_*$
since $t-cr<t_*$ for $t\in(t_*,t_*+cr)$.

Now we argue by contradiction, assume that $\phi_c(t)$ is not always decreasing
after $t_*$, then there exists a $t_1>t_*$ such that
$\phi_c(t)$ is strictly decreasing in $(t_*,t_1)$ and $\phi_c'(t_1)=0$.
Since $\phi_c(t)$ and $\phi_c^m(t)$ attain their local maximums at $t_*$
and local minimums at $t_1$, then
\begin{align*}
\phi_c'(t_*)=0, \quad (\phi_c^m)''(t_*)\le0, \quad d(\phi_c(t_*))\le b(\phi_c(t_*-cr)), \\
\phi_c'(t_1)=0, \quad (\phi_c^m)''(t_1)\ge0, \quad d(\phi_c(t_1))\ge b(\phi_c(t_1-cr)),
\end{align*}
and
$$
b(\phi_c(t_1-cr))\le d(\phi_c(t_1)) <d(\phi_c(t_*)) \le b(\phi_c(t_*-cr)),
$$
which implies that $\phi_c(t_1-cr)<\phi_c(t_*-cr)$ according to
the monotonicity of $b(\cdot)$.
Case I: if $t_1-t_*\le cr$, then $t_*-cr<t_1-cr\le t_*$,
and $\phi_c(t_1-cr)>\phi_c(t_*-cr)$ as $\phi_c(t)$ is strictly increasing before $t_*$.
We arrive at a contradiction.
Case II: if $t_1-t_*>cr$, then $t_1-cr> t_*$,
$\phi_c(t_1-cr)>\phi_c(t_1)$ as $\phi_c(t)$ is strictly decreasing in $(t_*,t_1)$,
and then
$$
b(\phi_c(t_1))<b(\phi_c(t_1-cr))\le d(\phi_c(t_1)).
$$
Noticing that $\phi_c(t_1)<\phi_c(t_*)<K$,
we have another contradiction to the fact
$b(s)>d(s)$ for all $s\in(0,K)$.
Therefore, $\phi_c(t)$ is always decreasing after $t_*$ before reaching zero,
and zero is the only equilibrium that smaller than $K$.
The proof is completed.
$\hfill\Box$

The monotonicity of birth function $b(\cdot)$ plays a crucial role in the proof
of the monotonicity of the sharp type local solution $\phi_c(t)$
in Lemma \ref{le-monotone}
and the monotonicity of $\phi_c(t)$ is the foundation of
the proof of the following monotone dependence of $\phi_c$ with respect to $c$.
Here we develop a phase transform approach to show
more precise behavior about $\phi_c(t)$.
We note that generally speaking this method is incapable of showing the existence
of traveling waves with time delay since the trajectories
with time delay may intersect with each other.
However, it provides a blueprint to draw more precise information
about the solutions if we have already known or assumed the existence.

\begin{lemma}[Monotone Dependence] \label{le-phic}
For any $c>0$, let $(0,t_*)$ be the maximal interval such that
$\phi_c(t)$ is strictly increasing and $\phi_c(t)<K$ on $(0,t_*)$
as defined in Lemma \ref{le-monotone}.
Then three cases can happen:\\ \indent
(i) $t_*<+\infty$ and $\phi_c(t_*)=K$, then $\phi_c'(t_*)\ge0$,
$\phi_c(t)$ may grow or decay after $t_*$;\\ \indent
(ii) $t_*<+\infty$ and $\phi_c(t_*)<K$, then
$\phi_c'(t_*)=0$ and $\phi_c(t)$ decays to zero and
never grows up to $K$ after $t_*$;\\ \indent
(iii) $t_*=+\infty$, then $\phi_c(t)$ is strictly increasing on $(0,+\infty)$
and $\phi_c(+\infty)\le K$, $\phi_c'(+\infty)=0$.\\
Hereafter, we only consider $\phi_c(t)$ within $(0,t_*)$.
Then $\phi_c(t)$ is strictly monotonically increasing with respect to $c$
in their joint interval.
\end{lemma}
{\it\bfseries Proof.}
The assertion about case (ii) follows from Lemma \ref{le-monotone}.
To prove the monotone dependence of $\phi_c$ with respect to $c$,
we develop a phase transform approach with time delay
inspired by \cite{HJMY,Jin-Yin}, see also \cite{JDE18XU,non-monotone}
for the proof of nonexistence of traveling waves by using this phase transform approach.
We note that Lemma \ref{le-semi-1} implies the locally monotone dependence
for $t$ near $0$,
but it is not enough as we need the globally monotone dependence.

From the construction of $\phi_c(t)$, we see that $\phi_c(t)$
is strictly increasing on $(0,t_*)$, either $\phi_c(t_*)=K$ or $\phi_c'(t_*)=0$,
and $\phi_c(t)$ satisfies the following equation
as we solve \eqref{eq-tw} locally step by step
\begin{equation} \label{eq-zlocal}
\begin{cases}
c\phi'(t)=D(\phi^m(t))''-d(\phi(t))+b(\phi(t-cr)), \quad t\in(0,t_*),\\
\phi(0)=0, \quad \phi(t) \text{~is strictly increasing on~} (0,t_*).
\end{cases}
\end{equation}
Let
$$\psi_c(t)=D(\phi_c^m(t))'=Dm\phi_c^{m-1}(t)\phi_c'(t).$$
Since $\phi_c(t)$ is strictly increasing on $(0,t_*)$,
\eqref{eq-zlocal} is transformed into
\begin{equation} \label{eq-zpp}
\begin{cases}
\displaystyle
\phi'(t)=\frac{\psi(t)}{Dm\phi^{m-t}(t)}, \\
\displaystyle
\psi'(t)=\frac{c\psi(t)}{Dm\phi^{m-t}(t)}
-\big(b(\hat\phi_{cr}(t))-d(\phi(t))\big),
\end{cases}
\end{equation}
where $\hat\phi_{cr}(t):=\phi_c(t-cr)$.

The key ingredient of the phase transform approach
is to write $t\in(0,t_*)$ as an inverse function
of $\phi\in(0,\phi_c(t_*))$ according to $\phi=\phi_c(t)$
and interpret $\phi_c(t-cr)$ as a function of $\phi$ depending on $\psi_c$.
That is, $\tilde t(\phi):=\phi_c^{-1}(\phi)$ for $\phi\in(0,\phi_c(t_*))$,
and $\tilde\psi(\phi):=\psi_c(\tilde t(\phi))$,
and further
\begin{equation} \label{eq-zphicr}
\tilde\phi_{cr}(\phi):=\hat\phi_{cr}(\tilde t(\phi))=
\phi_c(\tilde t(\phi)-cr)
=\inf_{\theta\ge0}\left\{\int_\theta^\phi\frac{Dms^{m-1}}{\tilde\psi(s)}ds\le cr\right\}.
\end{equation}
Since $\phi_c(t)$ is strictly increasing, we see that
$\tilde\phi_{cr}(\phi)<\phi$ according to \eqref{eq-zpp}
(or \eqref{eq-zphicr}).
The local solution $\phi_c(t)$ in $(0,t_*)$ corresponds to a trajectory
in the phase plane $(\phi,\psi)$ of \eqref{eq-zpp}
and $\tilde\psi(\phi)>0$ except for some possible isolated points $t_1$
such that $\phi_c'(t_1)=0$ and then $\tilde\psi(\phi_c(t_1))=0$.
Those points are isolated according to the strictly increasing property of $\phi_c(t)$.
As we have made change of the variables, we can write \eqref{eq-zpp} into
\begin{equation} \label{eq-zpsi}
\frac{d\tilde\psi}{d\phi}=c-\frac{Dm\phi^{m-1}(b(\tilde\phi_{cr}(\phi))-d(\phi))}{\tilde\psi},
\quad \phi\in(0,\phi_c(t_*)).
\end{equation}
From Lemma \ref{le-semi-1}, we know that $(\phi_c^m)'(0^+)=0$,
i.e., $\psi_c(0)=0$ and $\tilde\psi(0)=0$.

Let $\Gamma_c$ be the curve of
$$\bar\psi(\phi):=\frac{Dm\phi^{m-1}(b(\phi)-d(\phi))}{c}.$$
Then $\bar\psi(0)=\bar\psi(K)=0$, $\bar\psi(\phi)>0$ for $\phi\in(0,K)$,
$\bar\psi(\phi)$ is decreasing near $K$ and
$\Gamma_c$ divides $(0,K)\times (0,+\infty)$ into two parts,
$E_1:=\{(\phi,\psi);\phi\in(0,K), 0<\psi<\bar\psi(\phi)\}$
and $E_2:=((0,K)\times (0,+\infty))\backslash\overline E_1$.
Since $\tilde\phi_{cr}(\phi)<\phi$, we have
\begin{equation*}
\frac{d\tilde\psi}{d\phi}=c-\frac{Dm\phi^{m-1}(b(\tilde\phi_{cr}(\phi))-d(\phi))}{\tilde\psi}
>c-\frac{Dm\phi^{m-1}(b(\phi)-d(\phi))}{\tilde\psi}
\ge 0,
\end{equation*}
for $(\phi,\tilde\psi)\in E_2$.
In cases (ii) and (iii), $\psi_c(t_*)=0$ and then
$\tilde\psi(\phi_c(t_*))=0$.
That is, the trajectory lies below $\Gamma_c$ when $t$ nears $t_*$ and then
$$\sup_{\phi\in(0,\phi_c(t_*))}\tilde\psi(\phi)\le \sup_{\phi\in(0,K)}\bar\psi(\phi).$$
Otherwise, $\tilde\psi(\phi_c(t_*))$ will be positive, which contradicts to cases (ii) and (iii).

After settling down the generalized phase plane,
we now divide the proof of the monotonically increasing dependence of $\phi_c$
with respect to $c$ into two steps.
Let $c_1>c_2>0$ and $\phi_{c_1}(t)$, $\phi_{c_2}(t)$ be the local solutions
on the interval $(0,t_*(c_1))$ and $(0,t_*(c_2))$, respectively.
The symbols $t_*(c_1)$ and $t_*(c_2)$ are $t_*$ corresponding to $c_1$ and $c_2$
separately defined in Lemma \ref{le-monotone}.
The functions $\tilde t(\phi)$, $\tilde\psi(\phi)$ and $\tilde\phi_{cr}(\phi)$ defined above
corresponding to $c_1$ are denoted by $\tilde t_1(\phi)$,
$\tilde\psi_1(\phi)$, $\tilde\phi_{c_1r}(\phi;\tilde\psi_1)$.
Then the functions $\tilde t_2(\phi)$, $\tilde\psi_2(\phi)$, $\tilde\phi_{c_2r}(\phi;\tilde\psi_2)$
follow similarly corresponding to $c_2$.

Step I. We assert that $\phi_{c_1}(t)>\phi_{c_2}(t)$
for all $t\in(0,c_1r]$
and $\tilde\psi_1(\phi)>\tilde\psi_2(\phi)$ for $\phi\in(0,\phi_{c_1}(c_1r)]$.
For the phase plane corresponding to $c_1$, \eqref{eq-zpsi} reads
\begin{equation} \label{eq-zpsi-c1}
\frac{d\tilde\psi_1}{d\phi}=c_1-%
\frac{Dm\phi^{m-1}(b(\tilde\phi_{c_1r}(\phi;\tilde\psi_1))-d(\phi))}{\tilde\psi_1},
\quad \phi\in(0,\phi_{c_1}(t_*(c_1))),
\end{equation}
and the phase plane corresponding to $c_2$ follows similarly.
Within $(0,c_1r)$, we have
$\tilde\phi_{c_1r}(\phi;\tilde\psi_1)=\phi_{c_1}(\tilde t(\phi)-cr)=0$
and \eqref{eq-zpsi-c1} is actually
\begin{equation} \label{eq-zpsi-c1-1}
\frac{d\tilde\psi_1}{d\phi}=c_1+\frac{Dm\phi^{m-1}d(\phi)}{\tilde\psi_1},
\quad \phi\in(0,\phi_{c_1}(c_1r)).
\end{equation}
Therefore, for $\phi\in(0,\phi_{c_1}(c_1r))$ we have
$$
\frac{d\tilde\psi_1}{d\phi}=c_1+\frac{Dm\phi^{m-1}d(\phi)}{\tilde\psi_1}
>c_2+\frac{Dm\phi^{m-1}d(\phi)}{\tilde\psi_1}
$$
and
$$
\frac{d\tilde\psi_2}{d\phi}
=c_2-%
\frac{Dm\phi^{m-1}(b(\tilde\phi_{c_2r}(\phi;\tilde\psi_2))-d(\phi))}{\tilde\psi_2}
\le c_2+\frac{Dm\phi^{m-1}d(\phi)}{\tilde\psi_2}.
$$
The comparison principle for the above singular ODE implies
that $\tilde\psi_1(\phi)>\tilde\psi_2(\phi)$ for $\phi\in(0,\phi_{c_1}(c_1r))$,
see also \cite{HJMY,Jin-Yin} for details
by calculating ${d(\tilde\psi_1-\tilde\psi_2)}/{d\phi}$.

Here it should be noted that $\phi_{c_1}(t)>\phi_{c_2}(t)$ for all $t\in(0,c_1r]$
is not a simple conclusion of $\tilde\psi_1(\phi)>\tilde\psi_2(\phi)$ for $\phi\in(0,\phi_{c_1}(c_1r)]$
as the phase planes corresponding to $c_1$ and $c_2$ are different
(the changes of variables are different, i.e., $\tilde t_1(\phi)\not\equiv\tilde t_2(\phi)$).
We need to argue by contradiction.
Suppose that $\phi_{c_1}(t)>\phi_{c_2}(t)$ is not true for all $t\in(0,c_1r]$,
then there exists a $t_0\in(0,c_1r]$ such that $\phi_{c_1}(t_0)=\phi_{c_2}(t_0)$.
The choice of $t_0$ may not be unique, we choose the smallest one
as $\phi_{c_1}(t)>\phi_{c_2}(t)$ for $t$ near $0$ according to Lemma \ref{le-semi-1}.
Then $\phi_{c_1}'(t_0)\le\phi_{c_2}'(t_0)$,
$\phi_{c_1}(t_0)=\phi_{c_2}(t_0)=:\phi_0\in(0,\phi_{c_1}(c_1r)]$ and
$\tilde t_1(\phi_0)=t_0=\tilde t_2(\phi_0)$,
$$\tilde\psi_1(\phi_0)=m\phi_0^{m-1}\phi_{c_1}'(t_0)
\le m\phi_0^{m-1}\phi_{c_2}'(t_0)=\tilde\psi_2(\phi_0),$$
which contradicts to the fact that $\tilde\psi_1(\phi)>\tilde\psi_2(\phi)$ for all $\phi\in(0,\phi_{c_1}(c_1r)]$.

Step II.
We prove that $\phi_{c_1}(t)>\phi_{c_2}(t)$
for all $t\in(c_1r,\min\{t_*(c_1),t_*(c_2)\})$
and $\tilde\psi_1(\phi)>\tilde\psi_2(\phi)$ for
$\phi\in(\phi_{c_1}(c_1r),\min\{\phi_{c_1}(t_*(c_1)),\phi_{c_2}(t_*(c_2))\})$.
According to Step I, $\tilde\psi_1(\phi)>\tilde\psi_2(\phi)$ at $\phi=\phi_{c_1}(c_1r)$,
let $\phi_*>\phi_{c_1}(c_1r)$ be the first point such that
$\tilde\psi_1(\phi)>\tilde\psi_2(\phi)$ is not true
as we are arguing by contradiction.
Then $\tilde\psi_1(\phi)>\tilde\psi_2(\phi)$ for $\phi\in(0,\phi_*)$
and $\tilde\psi_1(\phi_*)=\tilde\psi_2(\phi_*)$,
$\tilde\psi_1'(\phi_*)\le\tilde\psi_2'(\phi_*)$.
For $\phi\in(\phi_{c_1}(c_1r),\phi_*)$, we have
$\phi>\phi_{c_1}(c_1r)>\phi_{c_2}(c_1r)$, $\tilde t_1(\phi)>c_1r>c_2r$,
$\tilde t_2(\phi)>c_1r>c_2r$,
and then \eqref{eq-zphicr} is simplified to
\begin{equation} \label{eq-zphi12}
\int_{\tilde\phi_{c_1r}(\phi;\tilde\psi_1)}^\phi
\frac{Dms^{m-1}}{\tilde\psi_1(s)}ds=c_1r,
\quad
\int_{\tilde\phi_{c_2r}(\phi;\tilde\psi_2)}^\phi
\frac{Dms^{m-1}}{\tilde\psi_2(s)}ds=c_2r.
\end{equation}
Since $c_1r>c_2r$ and $\tilde\psi_1(\phi)>\tilde\psi_2(\phi)$ for $\phi\in(0,\phi_*)$,
\eqref{eq-zphi12} tells us that
$$\tilde\phi_{c_1r}(\phi;\tilde\psi_1)<\tilde\phi_{c_2r}(\phi;\tilde\psi_2)$$
for $\phi\in(\phi_{c_1}(c_1r),\phi_*)$.

Now we use \eqref{eq-zpsi-c1} to deduce that
\begin{align*}
\frac{d\tilde\psi_1}{d\phi}=&c_1-%
\frac{Dm\phi^{m-1}(b(\tilde\phi_{c_1r}(\phi;\tilde\psi_1))-d(\phi))}{\tilde\psi_1}\\
>&c_1+\frac{Dm\phi^{m-1}d(\phi)}{\tilde\psi_1}-
\frac{Dm\phi^{m-1}b(\tilde\phi_{c_2r}(\phi;\tilde\psi_2))}{\tilde\psi_1}\\
>&c_1+\frac{Dm\phi^{m-1}d(\phi)}{\tilde\psi_1}-
\frac{Dm\phi^{m-1}b(\tilde\phi_{c_2r}(\phi;\tilde\psi_2))}{\tilde\psi_2},
\quad \phi\in(0,\phi_*),
\end{align*}
and similarly for the phase plane of $c_2$ we have
\begin{align*}
\frac{d\tilde\psi_2}{d\phi}=&c_2-%
\frac{Dm\phi^{m-1}(b(\tilde\phi_{c_2r}(\phi;\tilde\psi_2))-d(\phi))}{\tilde\psi_2}\\
=&c_2+\frac{Dm\phi^{m-1}d(\phi)}{\tilde\psi_2}-
\frac{Dm\phi^{m-1}b(\tilde\phi_{c_2r}(\phi;\tilde\psi_2))}{\tilde\psi_2},
\quad \phi\in(0,\phi_*).
\end{align*}
It follows that at the point $\phi_*$,
$$
\frac{d\tilde\psi_1}{d\phi}-\frac{d\tilde\psi_2}{d\phi}
>c_1-c_2>0,
$$
which is a contradiction to $\tilde\psi_1'(\phi_*)\le\tilde\psi_2'(\phi_*)$.
This argument by contradiction and Step I show that
$\tilde\psi_1(\phi)>\tilde\psi_2(\phi)$ for
$\phi\in(0,\min\{\phi_{c_1}(t_*(c_1)),\phi_{c_2}(t_*(c_2))\})$.
Using this fact, we can show that $\phi_{c_1}(t)>\phi_{c_2}(t)$
for all $t\in(0,\min\{t_*(c_1),t_*(c_2)\})$
in a similar procedure as in Step I.
The proof is completed.
$\hfill\Box$

\begin{lemma} \label{le-cstar}
There exists a unique $\hat c>0$ such that
$\phi_{\hat c}(t)$ strictly increasing on $(0,+\infty)$,
$\phi_{\hat c}(+\infty)=K$ and the function $\phi_{\hat c}(t)$ is also unique.
\end{lemma}
{\it\bfseries Proof.}
According to Lemma \ref{le-phic},
$\phi_c$ is strictly increasingly depending on $c>0$.
Lemma \ref{le-semi-decay} and Lemma \ref{le-semi-blow} show that
$\phi_c$ grows up to $K$ at finite time if $c$ is large
and decays to $0$ at finite time if $c$ is small.
Applying the continuous dependence Lemma \ref{le-continu}
and the monotone dependence Lemma \ref{le-phic} with respect to $c$,
we can define
\begin{equation} \label{eq-cstar}
\hat c=\{c>0; \phi_c(t) \text{~grows up to $K$ in finite time}\}.
\end{equation}
Then $\hat c>0$ according to Lemma \ref{le-semi-decay},
$\phi_{\hat c}$ is unique as we solve \eqref{eq-tw} step by step,
and $\hat c$ is the speed that satisfies the conditions in this Lemma.

We show that $\hat c$ is the unique speed that has the properties in this Lemma.
Lemma \ref{le-phic} implies the strictly monotone dependence of $\phi_c(t)$
with respect to $c$.
And more precisely, from the proof of Lemma \ref{le-phic}
we know that $\tilde\psi(\phi)$ is also strictly monotone dependence
with respect to $c$,
where $\tilde\psi(\phi)$ is the trajectory in the generalized phase plane
as in the proof of Lemma \ref{le-phic}.
The strictly monotone dependence implies the uniqueness of $\hat c$.
In fact, if there are $c_1>c_2$ that have the properties in this Lemma,
let $\tilde\psi_1(\phi)$ and $\tilde\psi_2(\phi)$
be the functions defined as in the proof of Lemma \ref{le-phic}.
Then $\tilde\psi_1(K)=\tilde\psi_2(K)$,
which contradicts to the strictly monotone dependence of
$\tilde\psi(\phi)$ with respect to $c$.
This completes the proof.
$\hfill\Box$

Lemma \ref{le-cstar} shows that $\phi_{\hat c}'(t)\ge0$ and
$\phi_{\hat c}(t)$ is strictly increasing in $(0,+\infty)$.
We need to prove a strong version as follows.

\begin{lemma} \label{le-phip}
The sharp traveling wave $\phi_{\hat c}(t)$ in Lemma \ref{le-cstar}
satisfies $\phi_{\hat c}'(t)>0$ for all $t\in(0,+\infty)$.
\end{lemma}
{\it\bfseries Proof.}
It is obvious that $\phi_{\hat c}'(t)\ge0$
and we argue by contradiction and assume that there exists a $t_0\in(0,+\infty)$
such that $\phi_{\hat c}'(t_0)=0$.
If $\phi_{\hat c}''(t_0)\ne0$, then $\phi_{\hat c}$ attains its
local strictly extreme value at $t_0$,
which cannot happen since $\phi_{\hat c}(t)$ is strictly increasing on $(0,+\infty)$.
Therefore, $\phi_{\hat c}''(t_0)=0$ and
\begin{equation} \label{eq-zexpand}
\phi_{\hat c}(t)=\phi_{\hat c}(t_0)+A(t-t_0)^3+o(|t-t_0|^3),
\quad t\to t_0,
\end{equation}
where $A\ge0$ as $\phi_{\hat c}(t)$ is strictly increasing.
The above expansion is valid as $\phi_{\hat c}$ is smooth away from the boundary of its support.
We first assume that $A>0$
and we write $\phi_{\hat c}(t)$ as $\phi(t)$ for simplicity.
Now \eqref{eq-tw} reads
\begin{align*}
b(\phi(t-cr))=&c\phi'(t)-D({\phi^m}(t))''+d(\phi(t))\\
=&3cA(t-t_0)^2-6DAm\phi^{m-1}(t_0)(t-t_0)+d(\phi(t_0))\\
&+3Ad'(\phi(t_0))(t-t_0)^2+o(|t-t_0|),
\quad t\to t_0,
\end{align*}
where the right hand side is monotonically decreasing near $t_0$.
By noticing that $\phi(t)$, $\phi(t-cr)$ and $b(\phi(t-cr))$
are strictly increasing, we arrive at a contradiction.
If $A=0$, we can expand \eqref{eq-zexpand} to higher odd order
and proceed the above argument similarly.
The proof is completed.
$\hfill\Box$

In order to compare $\hat c$ with $c^*(m,D,r)$ defined by \eqref{eq-def},
we need to compare the sharp type traveling wave $\phi_{\hat c}(t)$
with other smooth type traveling waves.

\begin{lemma} \label{le-compare}
There holds $\hat c=c^*(m,D,r)$.
That is, $\hat c$ is the minimal admissible traveling wave speed.
Furthermore, the wave speed for the smooth type traveling wave
is greater than the unique sharp type traveling wave speed $\hat c$.
\end{lemma}
{\it\bfseries Proof.}
Lemma \ref{le-cstar} shows that the sharp type traveling wave is unique.
We argue by contradiction in the following.
Let $\hat\phi(t)$ be a smooth type traveling wave with speed $c_1<\hat c$,
and let $\phi_{c_1}(t)$ and $\phi_{\hat c}(t)$
be the local sharp type solutions defined by \eqref{eq-semi}
corresponding to $c_1$ and $\hat c$ respectively.
Since $c_1<\hat c$ and $\phi_{\hat c}(t)$ increases in $(0,+\infty)$
with $\phi_{\hat c}(+\infty)=K$,
we see that $\phi_{c_1}(t_*(c_1))<K$ according to the strictly monotone dependence
and the uniqueness of $\hat c$ in Lemma \ref{le-cstar},
where $t_*(c_1)$ is defined in the proof of Lemma \ref{le-phic}.

Let $\tilde\psi(\phi)$ be the trajectory in the generalized phase plane
corresponding to the sharp type $\phi_{c_1}(t)$
defined in the proof of Lemma \ref{le-phic}.
Then $\tilde\psi(\phi_{c_1}(t_*(c_1)))=0$ as $\phi_{c_1}'(t_*(c_1))=0$.
For the monotonically increasing smooth type traveling wave solution $\hat\phi(t)$
we can also define the generalized phase plane
and let $\hat\psi(\phi)$ be the trajectory corresponding to the
smooth type traveling wave $\hat\phi(t)$.
The local asymptotic analysis Lemma \ref{le-semi-1} implies that
$$\tilde\psi(\phi)\sim c_1\phi, \quad \phi\to0^+,$$
and similar analysis shows that
$$\hat\psi(\phi)\sim \frac{Dm(b'(0)e^{-\lambda c_1r}-d'(0))}{c_1}\phi^m, \quad \phi\to0^+,$$
where $\lambda>0$ is the unique solution of the equation
$\lambda c_1+d'(0)=b'(0)e^{-\lambda c_1r}$.
The above local asymptotic behavior near zero shows that
$\hat\psi(\phi)<\tilde\psi(\phi)$ for $\phi\in(0,\phi_*)$ with some $\phi_*>0$.
Similar to the proof of Lemma \ref{le-phic},
we can show that
$\hat\psi(\phi)<\tilde\psi(\phi)$ for all $\phi\in(0,\phi_{c_1}(t_*(c_1))]$
and then $\hat\psi(\phi_0)=0$ for some $\phi_0\in(0,\phi_{c_1}(t_*(c_1)))\subset(0,K)$.
Similar to the proof of Lemma \ref{le-phip},
we can prove that $\hat\phi'(t)>0$ for all $t\in\mathbb R$
if $\hat\phi$ is a monotonically increasing smooth type traveling wave solution.
That is, $\hat\psi(\phi)>0$ for all $\phi\in(0,K)$,
which contracts to $\hat\psi(\phi_0)=0$ for some $\phi_0\in(0,K)$.

The above argument shows that the wave speed for the smooth type traveling wave
is greater than or equal to the unique sharp type traveling wave speed $\hat c$.
Next, we only need to show that there exist no smooth traveling waves
with speed $\hat c$.
The argument by contradiction is similar to the proof above
with the modifications such that $\phi_{c_1}(t_*(c_1))=K$
and $\hat\psi(\phi)<\tilde\psi(\phi)$ for all $\phi\in(0,K]$
and then $\hat\psi(\phi_0)=0$ for some $\phi_0\in(0,K)$.
This completes the proof.
$\hfill\Box$

The dependence of $c^*(m,D,r)$ with respect to the time delay $r$
is formulated via a variational characterization inspired by
Benguria and Depassier \cite{Benguria} and see also Huang et al. \cite{HJMY}.

\begin{lemma} \label{le-cstarr}
The minimal traveling wave speed $c^*(m,D,r)$ for the time delay $r>0$
is strictly smaller than that without time delay,
i.e., $c^*(m,D,r)<c^*(m,D,0)$.
\end{lemma}
{\it\bfseries Proof.}
Let $\phi(t)$ be the unique sharp type traveling wave
corresponding to the speed $c^*(m,D,r)=\hat c$ according to Lemma \ref{le-compare}.
In Lemma \ref{le-phic} we list three possible cases of the functions $\phi_c(t)$
corresponding to all $c>0$,
and combining the strictly monotone dependence (Lemma \ref{le-phic}) and
continuous dependence (Lemma \ref{le-continu}) of $\phi_c$
with respect to $c$ and the uniqueness of the sharp type traveling wave,
we see that $\phi(t)$ is a special function in case (iii)
and $\phi(t)$ is strictly increasing on $(0,+\infty)$,
$\phi(+\infty)=K$ and $\phi'(+\infty)=0$.
Lemma \ref{le-phip} shows that $\phi'(t)>0$ for all $t\in(0,+\infty)$.

In the proof of Lemma \ref{le-phic},
we develop the generalized phase plane \eqref{eq-zpp} and \eqref{eq-zpsi}
as $\tilde\phi_{cr}(\phi)$ is defined by \eqref{eq-zphicr}.
And additionally, $\tilde\psi(0)=0$, $\tilde\psi(K)=0$ since $\phi'(+\infty)=0$,
$\tilde\psi(\phi)>0$ for all $\phi\in(0,K)$
since $\phi'(t)>0$ for all $t\in(0,\infty)$.
Now, we rewrite \eqref{eq-zpsi} into
\begin{equation} \label{eq-zpsi-vc}
\frac{d\tilde\psi}{d\phi}=c
-\frac{Dm\phi^{m-1}(b(\phi)-d(\phi))}{\tilde\psi}
+\frac{Dm\phi^{m-1}(b(\phi)-b(\tilde\phi_{cr}(\phi)))}{\tilde\psi},
~\phi\in(0,K).
\end{equation}
For any $g\in\mathscr{D}=\{g\in C^1([0,K]);g(K)=0,\int_0^K g(s)ds=1,g'(s)<0,\forall s\in(0,K)\}$,
we multiply \eqref{eq-zpsi-vc} by $g(s)$ and integrate over $(0,K)$ to find
\begin{align} \nonumber
c=&\int_0^K g(\phi)\frac{d\tilde\psi}{d\phi}d\phi
+\int_0^K g(\phi)\frac{Dm\phi^{m-1}(b(\phi)-d(\phi))}{\tilde\psi}d\phi
\\ \nonumber
&-\int_0^K g(\phi)\frac{Dm\phi^{m-1}(b(\phi)-b(\tilde\phi_{cr}(\phi)))}{\tilde\psi}d\phi
\\ \nonumber
=&\int_0^K -g'(\phi)\tilde\psi(\phi)d\phi
+\int_0^K g(\phi)\frac{Dm\phi^{m-1}(b(\phi)-d(\phi))}{\tilde\psi}d\phi
\\ \nonumber
&+\big[g(\phi)\tilde\psi(\phi)\big]\Big|_{\phi=0}^{\phi=K}
-\int_0^K g(\phi)\frac{Dm\phi^{m-1}(b(\phi)-b(\tilde\phi_{cr}(\phi)))}{\tilde\psi}d\phi
\\ \nonumber
\ge&
2\sqrt{D}\int_0^K \sqrt{-ms^{m-1}g(s)g'(s)(b(s)-d(s))}ds
\\ \label{eq-zvar}
&-\int_0^K g(\phi)\frac{Dm\phi^{m-1}(b(\phi)-b(\tilde\phi_{cr}(\phi)))}{\tilde\psi}d\phi,
\end{align}
by Cauchy inequality and according to
$\int_0^K g(s)ds=1$ and $[g(\phi)\tilde\psi(\phi)]\big|_{\phi=0}^{\phi=K}=0$ as
$\tilde\psi(0)=0=g(K)$.

It should be noted that the equality in \eqref{eq-zvar} is attainable
at some function $\hat g$ such that
\begin{equation} \label{eq-zhatg}
-\hat g'(\phi)\tilde\psi(\phi)=\hat g(\phi)\frac{Dm\phi^{m-1}(b(\phi)-d(\phi))}{\tilde\psi},
\quad \phi\in(0,K),
\end{equation}
with $\hat g(K)=0$ and $\hat g'(\phi)<0$ for all $\phi\in(0,K)$.
In fact, such kind of solution $\hat g$ to \eqref{eq-zhatg} is solvable
since $\tilde\psi(K)=0$ and $\tilde\psi(\phi)\sim \kappa(K-\phi)$ as $\phi\to K^-$
for some $\kappa>0$
and $\tilde\psi(\phi)>0$ for all $\phi\in(0,K)$ according to Lemma \ref{le-phip},
see for example \cite{HJMY} for the phase plane without time delay.
Then
$$
(\ln \hat g)'=\frac{\hat g'}{\hat g}=-
\frac{Dm\phi^{m-1}(b(\phi)-d(\phi))}{\tilde\psi^2(\phi)}
\sim
-\frac{DmK^{m-1}(d'(K)-b'(K))}{\kappa^2(K-\phi)},
~\phi\to K^-,
$$
which has infinitely many solutions with $\hat g(K)=0$.
(Otherwise, if $\tilde\psi(K)\ne0$, then $\hat g(\phi)\equiv0$ is the unique solution
of \eqref{eq-zhatg} with $\hat g(K)=0$.)
On the other hand, in the proof of Lemma \ref{le-compare},
we show that
$\tilde\psi(\phi)\sim c\phi$ as $\phi\to0^+$,
and hence
$$
(\ln \hat g)'=\frac{\hat g'}{\hat g}=-
\frac{Dm\phi^{m-1}(b(\phi)-d(\phi))}{\tilde\psi^2(\phi)}
\sim
-\frac{Dm(b'(0)-d'(0))\phi^{m-2}}{c^2},
~\phi\to 0^+.
$$
It follows that $\hat g(0)<+\infty$ as $m>1$ and $\hat g\in \mathscr D$
such that the equality in \eqref{eq-zvar} is attainable.
(Here we point out that for the smooth type traveling waves,
$\hat\psi(\phi)\sim \mu \phi^m$ as $\phi\to0^+$ for some $\mu>0$
as in the proof of Lemma \ref{le-compare},
and then $\hat g(0)=+\infty$ in this case and
the equality in \eqref{eq-zhatg} is not attainable.)

Now that we showed the attainable of the equality in \eqref{eq-zvar} at $\hat g$,
we have
\begin{align} \nonumber
c=&2\sqrt{D}\int_0^K \sqrt{-ms^{m-1}\hat g(s)\hat g'(s)(b(s)-d(s))}ds
\\ \nonumber
&-\int_0^K \hat g(\phi)\frac{Dm\phi^{m-1}(b(\phi)-b(\tilde\phi_{cr}(\phi)))}{\tilde\psi}d\phi
\\ \nonumber
<& 2\sqrt{D}\int_0^K \sqrt{-ms^{m-1}\hat g(s)\hat g'(s)(b(s)-d(s))}ds
\\ \label{eq-zc}
\le& \sup_{g\in\mathscr D}2\sqrt{D}\int_0^K \sqrt{-ms^{m-1}\hat g(s)\hat g'(s)(b(s)-d(s))}ds,
\end{align}
where the ``$<$'' in \eqref{eq-zc} follows from the
strictly monotone increasing of $\phi(c)$ such that $\tilde\phi_{cr}<\phi$
as in the proof of Lemma \ref{le-phic}.
The proof is completed.
$\hfill\Box$

{\it\bfseries Proof of Theorem \ref{th-existence}.}
The uniqueness of the wave speed of the sharp type traveling wave
is proved in Lemma \ref{le-cstar}.
We see that the corresponding sharp wave is also unique (up to shift)
since the maximal solution $\phi_c^1(t)$ in Lemma \ref{le-semi-1} is the unique solution
such that $\phi_c^1(t)>0$ in a right neighbor of $0$.
The monotonicity follows from Lemma \ref{le-monotone}, Lemma \ref{le-phic} and
Lemma \ref{le-cstar}.
Lemma \ref{le-compare} implies that $c^*(m,D,r)=\hat c$.
The positiveness of $c^*(m,D,r)$ and the uniqueness of $c^*(m,D,r)$ and the
sharp type traveling wave follow from Lemma \ref{le-cstar}.
The dependence of $c^*(m,D,r)$ with respect to $r$
is proved in Lemma \ref{le-cstarr}.
The proof is completed by combining these lemmas.
$\hfill\Box$

{\it\bfseries Proof of Theorem \ref{th-smooth}.}
This is proved in Lemma \ref{le-compare}.
The regularity is trivial since $\phi(t)>0$ for all $t\in\mathbb R$,
where \eqref{eq-tw} is non-degenerate.
$\hfill\Box$

{\it\bfseries Proof of Theorem \ref{th-sharp}.}
The asymptotic behavior near $0$ in Lemma \ref{le-semi-1}
completes the proof,
see also \cite{non-monotone} for details.
$\hfill\Box$

\section*{Acknowledgement}

This work was done when T.Y. Xu and S.M. Ji visited McGill University
supported by CSC programs.
They would like to express their sincere thanks for the hospitality
of McGill University and CSC.
The research of S. Ji was
supported by NSFC Grant No.~11701184 and CSC No. 201906155021
and the Fundamental Research Funds for the Central Universities of SCUT.
The research of M. Mei was supported in part
by NSERC Grant RGPIN 354724-16, and FRQNT Grant No. 2019-CO-256440.
The research of J. Yin was supported in
part by NSFC Grant No. 11771156 and NSF of Guangzhou Grant No. 201804010391.

\end{document}